\documentclass[12pt]{article}
\usepackage{amsmath,amsthm,amssymb, amstext, amsopn, amsxtra}
\usepackage{latexsym}
\usepackage{amsfonts}
\begin{document}
\newcommand{\nt}{\noindent}
\newcommand{\bs}{\bigskip}
\newcommand{\ms}{\medskip}
\newcommand{\mk}{\medskip}
\newcommand{\sk}{\smallskip}
\newcommand{\ep}{\varepsilon}
\newcommand{\m}{{\mathfrak m}}
\newcommand{\p}{{\mathfrak p}}
\newcommand{\dd}{\delta}
\newcommand{\A}{{\mathbb A}}
\newcommand{\R}{{\mathcal R}}
\newcommand{\M}{{\mathbb M}}
\newcommand{\rop}{\rho^{\,\prime}(R)}
\newcommand{\pq}{\langle\,p,q\,\rangle}
\newcommand{\LA}{\langle\!\langle\,}
\newcommand{\RA}{\,\rangle\!\rangle}
\newcommand{\QQ}{[\,2\,]}
\[   \]

\begin{center}
\Large{\bf A Quantum-Trace Determinantal Formula for\\
Matrix Commutators, and Applications}

\bigskip
\large{Dinesh Khurana, T.\,Y.~Lam, and Noam Shomron} 
\end{center}

\bigskip

\begin{abstract}
\begin{small}
In this paper, we establish a determinantal formula for $\,2\times 2\,$
matrix commutators $\,[X,Y]=XY-YX\,$ over a commutative ring, using 
(among other invariants) the quantum traces of $\,X\,$ and $\,Y$.  
Special forms of this determinantal formula include a ``trace version'', 
and a ``supertrace version''.  Some applications of these formulas are
given to the study of value sets of binary quadratic forms, the 
factorization of $\,2\times 2\,$ integral matrices, and the solution 
of certain simultaneous diophantine equations over commutative rings.
\end{small}
\end{abstract}

\sk

\bs\nt 
{\bf \S1. \ Introduction}

\bs
Throughout this paper, the word ``commutator'' is taken to mean an
{\it additive\/} commutator $\,[X,Y]:=XY-YX$, where $\,X,Y\,$ are 
elements in an (associative) ring $\,R$.  The case of special interest 
to us in this paper is when $\,R\,$ is a matrix ring $\,{\mathbb M}_n(S)$, 
where $\,S\,$ is a commutative ring (with identity).

\bs
In the case $\,S\,$ is a field (or a division ring), several papers 
in the literature have dealt with the theme of computing the possible 
ranks of $\,[X,Y]$, when $\,X\,$ is a given matrix. An important  
special case of this is to determine when $\,[X,Y]\,$ can achieve full 
rank for a suitable $\,Y$; see, e.g.~[GL] and [S\'a].  More recently, 
two of the present authors have taken up the study of similar themes 
when the division ring $\,S\,$ is replaced by a more general ring.  
Here, it is no longer possible to use effectively the notion of matrix 
ranks.  But the classical full-rank case corresponds to the {\it 
invertibility\/} of the commutator $\,[X,Y]$, which is certainly a 
condition of interest and significance; see [KL$_1$, KL$_2$].

\bs
Over a commutative base ring $\,S$, a matrix $\,M\in {\mathbb M}_n(S)\,$
is invertible iff $\,{\rm det}\,(M)\,$ is a unit in $\,S$. Thus,
deciding the invertibility of $\,[X,Y]\,$ rests on understanding the 
behavior of its determinant $\,{\rm det}\,[X,Y]$.  However, there 
seems to be no general formula available in the literature for the 
computation of such an $\,n\times n\,$ determinant. It is much more 
feasible to find nice formulas in the $2\times 2$ case, since 
$2\times 2$ determinants are so much easier to compute.  Indeed, by 
using the Cayley-Hamilton theorem, or using classical adjoints, several 
$2\times 2$ formulas can be written down for $\,{\rm det}\,[X,Y]$.  
For a quick survey of this, see \S2.  These formulas often involve 
higher powers of $\,X\,$ and $\,Y$; unfortunately, this feature tends 
to greatly limit their applicability, for instance, to the study of 
the invertibility of $\,[X,Y]$. Optimally, we should hope to write 
down formulas expressing $\,{\rm det}\,[X,Y]\,$ in terms of quantities 
(e.g.~traces and determinants) directly associated with the matrices 
$\,X,Y,XY$, and if necessary, $\,YX$.

\bs
In the beginning phase of our work, we were aware of the following 
``trace version'' of a $\,2\times 2\,$ determinantal formula for
$\,{\rm det}\,[X,Y]\,$ in the context of P.I.-theory and invariant theory: 
if $\,X,Y\in {\mathbb M}_2(S)\,$ have traces $\,t,\,t^{\,\prime}\,$ and 
determinants $\,\dd,\,\dd^{\,\prime}$, then
$$
{\rm det}\;[X,Y] = 4\,\dd^{\,\prime}\dd - (\,{\rm tr}\,(XY)\,)^2
- \dd\,t^{\,\prime\,2} - \dd^{\,\prime}\,t^2 
+ {\rm tr}\,(XY)\;t^{\,\prime}t\,.  \leqno (1.1)
$$
While this formula is a natural consequence of the invariant theory of
$\,2\times 2\,$ matrices (under simultaneous conjugation), it turned out 
to be not the most suitable for the applications we have in mind.  In 
the search for a better alternative (for our purposes), we stumbled 
upon a ``super-trace version'' of a determinantal formula, using the 
``supertrace'' of $2\times 2$ matrices $\,M=(m_{ij})$, which is defined 
to be $\,{\rm str}\,(M)=m_{11}-m_{22}\in S$. (The supertrace terminology 
and its concomitant notation ``str'' come from the theory of 
super-algebras.)  If $\,X,Y\in {\mathbb M}_2(S)\,$ have determinants 
$\,\dd,\,\dd^{\,\prime}$ and supertraces $\,\tau,\,\tau^{\,\prime}$, 
the ``supertrace version'' of the commutator determinantal formula 
states the following:
$$
-\,{\rm det}\;[X,Y] = \dd\,\tau^{\,\prime\,2} + \dd^{\,\prime}\,\tau^2
+{\rm tr}\,(XY)\;\tau^{\,\prime}\tau-{\rm str}\,(XY)\;{\rm str}\,(YX).
\leqno (1.2)
$$
Contrary to the case of $(1.1)$, the existence of $(1.2)$ does not 
seem to be predictable by invariant theory since the supertrace of
a matrix in $\,{\mathbb M}_2(S)\,$ is not a similarity invariant.
Thus, in a manner of speaking, the existence of the formula (1.2) is
a bit surprising.

\bs
While the two formulas (1.1) and (1.2) look substantially different 
(having for instance a different number of terms), they do share some 
common features. This begs the question whether they are special cases 
of one single more general formula. In the world of quantum mathematics, 
it is easy to speculate what might be the case.  If $\,q\,$ is a fixed 
element in the ring $\,S$, there is a well-known notion of a 
$\,q$-{\it trace\/} for $\,M=(m_{ij})\in {\mathbb M}_2(S)$, defined by 
$\,{\rm tr}_q(M)=m_{11}+q\,m_{22}$. For $\,q={\pm 1}$, this $\,q$-trace 
retrieves the trace and the supertrace respectively. The existence 
of (1.1) and (1.2) would seem to strongly suggest that there is a hybrid
``quantum version'' of a determinantal formula for $\,{\rm det}\,[X,Y]$, 
which would specialize to (1.1) when $\,q=1$, and to (1.2) when $\,q=-1$.

\bs
In \S4 of this paper, we prove that this is indeed the case. The main
result is given in Theorem 4.1 which, for any given $\,q$, expresses
$\,q\cdot {\rm det}\,[X,Y]\,$ in terms of various quantities, including
the $\,q$-traces of $\,X,Y,XY$, and $\,YX$. This theorem is preceded by 
Theorem 3.4 in \S3, which treats the special case of $\,q$-traceless 
matrices $\,X,Y$.

\bs
Some applications of the determinantal formulas are given in Sections
5--6.  In \S5, we present a new characteristic-free treatment of a norm 
theorem of Taussky [Ta$_1$, Ta$_2$] for quadratic field extensions, and 
give an explicit generic version of this result for commutative rings.  
This is followed by a last section (\S6) devoted to the further study of 
the value sets of binary quadratic forms, and factorization questions on 
$2\times 2$ matrices over rings.  For given elements $\,p,q\,$ in a 
commutative ring $\,S$, we show in Theorem 6.3 that 
{\it a non $0$-divisor $\,c\in S\,$ can be written in the form 
$\,pr^2+qs^2\,$ for some $\,r,s\in S\,$ iff the matrix $\begin{pmatrix}
0&c\,q\\-c\,p&0\end{pmatrix}$ can be factored into $\,XY\,$ such that 
$\,{\rm det}\,(X)=c\,p$, $\,{\rm det}\,(Y)=c\,q$, and $\,{\rm det}\,[X,Y]
=-c^2$.} Using this factorization theorem, we prove a result (Thm.~6.8) 
on affine curves over rings which implies that, for any element 
$\,c\in S\,$ representable in the form $\,pr^2+qs^2$, there exist 
$\,x,y,z\in S\,$ such that $\,px+qy=-c\,$ and $\,xy-z^2=-c^2$. Further 
applications of the determinantal formulas to the study of invertible 
commutators in matrix rings can be found in the forthcoming work 
[KL$_1$].

\mk

\bs\nt
{\bf \S2. \ Preliminary Determinantal Formulas}

\bs
For the rest of this paper, $R\,$ denotes the matrix ring 
$\,{\mathbb M}_2(S)$, where $\,S\,$ is a fixed commutative ring.  
In this section, we put together a few general (mostly known) facts 
pertaining to the computation of the determinant of a $2\times 2$ matrix 
commutator $\,[X,Y]=XY-YX$, where $\,X,Y\in R$.  This quick survey will 
pave the way to our more detailed investigations on 
$\,{\rm det}\,[X,Y]\,$ in the ensuing sections.

\bs
To begin with, we note that, since $\,[X,Y]\,$ is a $2\times 2$ traceless 
matrix for all $\,X,Y\in R$, the Cayley-Hamilton Theorem implies that
$$
[X,Y]^2=-\,{\rm det}\,[X,Y]\cdot I_2.\leqno (2.1)
$$
Thus, computing $\,{\rm det}\,[X,Y]\,$ is tantamount to computing the 
scalar matrix $\,-[X,Y]^2$.  This observation leads to the first 
available formula for $\,{\rm det}\,[X,Y]$.

\bs\nt
{\bf Proposition 2.2.} {\it For any $\,X,Y\in R$,}
$\,{\rm det}\;[X,Y]={\rm tr}\,\bigl(X^2Y^2\bigr)
-{\rm tr}\,\bigl( (XY)^2 \bigr)$.

\bs\nt
{\bf Proof.} We may assume that $\,S\,$ is the free commutative 
$\,{\mathbb Z}$-algebra generated by the eight entries of $\,X\,$ and
$\,Y$. Computing the trace on the matrices in (2.1), we have
$$
2\;{\rm det}\,[X,Y] = -\,{\rm tr}\,\bigl( [X,Y]^2\bigr)
      = -\,{\rm tr}\,\bigl(XYXY+YXYX-XYYX-YXXY\bigr).\leqno (2.3)
$$
On the RHS, $\,X\,(YXY)\,$ and $\,(YXY)\,X\,$ have the same trace, and
$\,(XY)\,(YX)\,$ and $\,(YX)\,(XY)\,$ have the same trace. Therefore,
the RHS of (2.3) can be rewritten as
$$
-2\;{\rm tr}\,\bigl( (XY)^2 \bigr)+2\;{\rm tr}\,(XYYX)
= -2\;{\rm tr}\,\bigl( (XY)^2 \bigr)+2\;{\rm tr}\,\bigl(X^2Y^2\bigr).
$$
Cancelling the factors of $\,2\,$ gives the desired result.\qed

\bs
Next, we'll prove a determinantal formula involving the {\it adjoint 
matrix.}  For a $2\times 2$ matrix $\,M=\begin{pmatrix}a&b\\c&d
\end{pmatrix}\in R$, we'll write $\,M^{\,\prime}=\begin{pmatrix}d&-b\\
-c&a\end{pmatrix}$ for the classical adjoint of $\,M$.  For 
$\,E=\begin{pmatrix} 0&-1\\1&0\end{pmatrix}$, we have the relation 
$\,M^{\,\prime}=EM^TE^{-1}$, as was noted in [GK:~(2.4)(3)]. This 
implies, in particular, that $\,(MN)^{\,\prime}
=N^{\,\prime}M^{\,\prime}\,$ for $\,M,N\in R$. 
The following lemma is a useful observation on the adjoint.

\bs\nt
{\bf Lemma 2.4.} {\it For any matrices $\,A,B\in R={\mathbb M}_2(S)$, 
we have}
$$
{\rm det}\;(A-B)={\rm det}\,(A)+{\rm det}\,(B)-{\rm tr}\,(AB^{\,\prime}).
\leqno (2.5)
$$

\sk\nt
{\bf Proof.} For any $\,N\in R$, we have $\,N\cdot N^{\,\prime}
=({\rm det}\;N)\cdot I_2$, and $\,N+N^{\,\prime}=({\rm tr}\;N)\cdot I_2$.
Applying these facts to the equations
\begin{eqnarray*}
(A-B)\,(A-B)^{\,\prime} \!&=&\! A\,A^{\,\prime}+B\,B^{\,\prime}
-A\,B^{\,\prime}-B\,A^{\,\prime}\\
\!&=&\! A\,A^{\,\prime}+B\,B^{\,\prime}-\bigl( A\,B^{\,\prime}
+ (A\,B^{\,\prime})^{\,\prime}\,\bigr),
\end{eqnarray*}
we arrive at the desired equation.\qed

\bs\nt
{\bf Proposition 2.6.} {\it For any matrices $\,X,Y\in R$, we have}
$$
{\rm det}\;[X,Y] = 2\;({\rm det}\;X)\,({\rm det}\;Y)
-{\rm tr}\;(XYX^{\,\prime}Y^{\,\prime}).\leqno (2.7)
$$

\sk\nt
{\bf Proof.} In the formula (2.5), let $\,A=XY\,$ and $\,B=YX$. Then
the LHS of (2.5) is $\,{\rm det}\;[X,Y]$, and its RHS is
$$
\quad
{\rm det}\,(XY)+{\rm det}\,(YX)-{\rm tr}\,\bigl( (XY)\,(YX)^{\,\prime}\bigr)
= 2\;({\rm det}\,X)\,({\rm det}\,Y)-{\rm tr}\,(XYX^{\,\prime}Y^{\,\prime}).
\quad  \qed
$$

\sk
The formulas in (2.2) and (2.7) are interesting, but are not very suitable 
for practical computations since their RHS's involve new matrices such 
as $\,(XY)^2$,  $\,X^2Y^2$, and $\,XYX^{\,\prime} Y^{\,\prime}$. Ideally, 
we would like to have formulas that express $\,{\rm det}\;[X,Y]\,$ in 
terms of quantities naturally associated with the matrices $\,X,Y,XY$, 
and $\,YX$. Such formulas will be obtained in the next two sections.
Nevertheless, combining the two determinantal formulas obtained so far
leads to the following curious trace identity.

\bs\nt
{\bf Corollary 2.8.} {\it For $\,X,Y\in R$,}
$\,\bigl( {\rm tr}\,(XY)\bigr)^2={\rm tr}\,(X^2Y^2)
+{\rm tr}\,\bigl(XYX^{\,\prime}Y^{\,\prime}\bigr)$.

\bs\nt
{\bf Proof.} Equating the two expressions for $\,{\rm det}\,[X,Y]\,$ in
(2.2) and (2.7), we get
$$
{\rm tr}\,\bigl(X^2Y^2\bigr)
-{\rm tr}\,\bigl( (XY)^2 \bigr)
= 2\;({\rm det}\;X)\,({\rm det}\;Y) -{\rm tr}\;(XYX^{\prime}Y^{\prime}).
\leqno (2.9)
$$
Using the well-known identity $\,2\;{\rm det}\,(M)
=\bigl({\rm tr}\,(M)\bigr)^2-{\rm tr}\,(M^2)\,$ for $\,M=XY\,$ on the 
RHS, we can cancel the terms $\,-{\rm tr}\,\bigl( (XY)^2\bigr)\,$ from 
both sides of (2.9). After this, transposition yields the desired 
result.\qed

\mk

\bs\nt
{\bf \S3. \ Determinantal Formula: The \,q\,-Traceless Case}

\bs
In this section, we shall embark upon the task of finding a determinantal
formula for $\,{\rm det}\,[X,Y]\,$ that is suitable for the applications 
we have in mind.  At the beginning stage of our work, we had at our
disposal both a ``trace version'' and a ``supertrace version'' of such a 
formula. This suggested to us that there is perhaps a ``$q$-trace version'' 
of the formula, which, for $\,q=1\,$ and $\,q=-1\,$ respectively, would
specialize to the trace version and the supertrace version.  Such a
quantum-trace version would then provide some kind of a ``homotopy''
from the trace version to the supertrace version, and conversely. 
In this and the following section, we shall begin the work to derive 
such a $\,q$-trace determinantal formula for $\,{\rm det}\,[X,Y]$. 
{\it Throughout these two sections, the element $\,q\in S\,$ will be 
regarded as fixed.}  Recall that, for any matrix $\,M=(m_{ij})\in 
{\mathbb M}_2(S)$, the $\,q$-trace $\,{\rm tr}_q(M)\,$ is defined to be 
$\,m_{11}+q\,m_{22}\in S$.

\bs
The strategy of our approach is to first find a $\,q$-trace determinantal
formula {\it in the case where $\,X,\,Y\,$ are both $\,q$-traceless\/}; 
that is, where $\,{\rm tr}_q(X)={\rm tr}_q(Y)=0$.  To guess how such a 
formula might look like, it would be a good idea to start with the case 
of the ordinary trace; that is, where $\,q=1$.  Here, a formula for 
$\,{\rm det}\,[X,Y]\,$ can be derived from the one in Proposition 2.2;
see, e.g.~the proof of [BBO:~Lemma 6].

\bs\nt
{\bf Proposition 3.1.} {\it Let $\,X,Y\in R={\mathbb M}_2(S)\,$ be such 
that $\,{\rm tr}\,(X)={\rm tr}\,(Y)=0$. Then}
$$
{\rm det}\,\;[X,Y]=4\;{\rm det}\,(XY)-\bigl( {\rm tr}\,(XY)\bigr)^2.
\leqno (3.2)
$$

\mk\nt
{\bf Proof.} It is of interest to observe that the RHS of (3.2) is just 
the negative of the usual discriminant of the quadratic characteristic 
polynomial of the product matrix $\,XY$. Let $\,\dd={\rm det}\,(X)\,$ and
$\,\dd^{\,\prime}={\rm det}\,(Y)\,$ (so that $\,{\rm det}\,(XY)
=\dd^{\,\prime}\dd\,$).  Since $\,X,Y\,$ are traceless, Cayley-Hamilton 
gives $\,X^2=-\dd\,I_2\,$ and $\,Y^2=-\dd^{\,\prime}\,I_2$. 
Thus, using (2.2) and the identity for $\,2\;{\rm det}\,(M)\,$ in the 
proof of (2.8), we have 
$$
{\rm det}\,[X,Y]={\rm tr}\,\bigl(\dd^{\,\prime}\dd\,I_2 \bigr)
-{\rm tr}\,\bigl((XY)^2\bigr)= 2\,\dd^{\,\prime}\dd
- \bigl( {\rm tr}\,(XY)\bigr)^2+ 2\;{\rm det}\,(XY),
$$
which simplifies to the RHS of (3.2).\qed

\bs\nt
{\bf Remark 3.3.} We could have also proved (3.2) by applying the methods 
of Procesi in [Pr$_1$]. From Procesi's approach (see [KP:~\S2.4, p.\,21]),
it will be enough to prove (3.2) in the case where the first traceless 
matrix $\,X\,$ is {\it diagonal\/}.  In this case, (3.2) can be easily 
checked by a direct computation of both sides of the equation.

\bs
To generalize Proposition 3.1 to the case of $\,q$-traces, we'll use the
following standard notation in quantum computations: with $\,q\in S\,$
fixed, we write 
$$
[\,n\,]=[\,n\,]_q=(q^n-1)/(q-1)=q^{n-1}+\cdots+q+1\in S\;\;\;(\,\forall\,
n\in {\mathbb Z}^{+}).
$$
In fact, we'll use this notation only for $\,n=2$; namely, $\,\QQ=1+q$.  
Treating this as a ``quantum $2$'', we could use $\,\QQ^2\,$ to replace
the factor $\,4\,$ in the formula (3.2).  Also, for $\,q$-traces, 
$\,{\rm tr}_q(XY)\,$ and $\,{\rm tr}_q(YX)\,$ may no longer be the same, 
so it would be wise to replace the term $\,\bigl( {\rm tr}\,(XY)\bigr)^2\,$ 
by $\,{\rm tr}_q(XY)\,{\rm tr}_q(YX)$. Fortuitously, these ``replacements'' 
turn out to give the following correct generalization of Proposition 3.1. 

\bs\nt
{\bf Theorem 3.4.} {\it Let $\,X,Y\in R={\mathbb M}_2(S)\,$ be such that 
$\,{\rm tr}_q(X)={\rm tr}_q(Y)=0$. Then}
$$
q\cdot {\rm det}\,[X,Y]=\QQ^2\,{\rm det}\,(XY)
-{\rm tr}_q(XY)\, {\rm tr}_q(YX).\leqno (3.5)
$$
{\it In this formula, the first term on the\/} RHS {\it can also be 
replaced by $\,{\rm tr}_q(X^2)\,{\rm tr}_q(Y^2)$.}

\bs\nt
{\bf Proof.} To prove the last statement, we note again that
Cayley-Hamilton implies $\,X^2=t\,X-{\rm det}\,(X)\,I_2$. Thus, taking 
$\,q$-traces gives $\,{\rm tr}_q(X^2)=-\QQ \cdot {\rm det}\,(X)$.  
Multiplying this by a similar equation for $\,Y$, we see that
$\,{\rm tr}_q(X^2)\,{\rm tr}_q(Y^2)=\QQ^2\,{\rm det}\,(XY)$.

\mk
To prove (3.5), let $\,X=\begin{pmatrix} -qd&b\\c&d\end{pmatrix}$, and
$\,Y=\begin{pmatrix}-qh&f\\g&h\end{pmatrix}$.  By direct computation, 
$\,[X,Y]\,$ has the form $\begin{pmatrix} k&\QQ\,r\\ \QQ\,r^{\,\prime}&-k 
\end{pmatrix}$, where $\,k=bg-cf$, $\,r=bh-df$, and $\,r^{\,\prime}
=dg-ch$.  Therefore, the LHS of (3.5) is
$$
q\cdot {\rm det}\,[X,Y]=-q\cdot \bigl(\QQ^2\,r^{\,\prime}r + k^2\bigr).
\leqno (3.6)
$$
To compute the RHS of (3.5), note that, in terms of descending powers of 
$\,\QQ\,$:
$$
{\rm tr}_q(XY)=q^2dh+q\,(cf+dh)+bg=\QQ^2\,dh+\QQ\,s+k, \;\,
\mbox{where}\;\;s=cf-dh.\leqno (3.7)
$$   
Similarly, we have 
$\,{\rm tr}_q(YX)=\QQ^2\,dh+\QQ\,s^{\,\prime}-k$, where $\,s^{\,\prime}
=bg-dh$. (For later use, note that $\,s^{\,\prime}-s=k$.) On the other
hand, 
$$
{\rm det}\,(XY)=(qd^2+bc)\,(qh^2+fg)
= \bigl(\,\QQ\,d^2+bc-d^2\bigr)\, \bigl(\,\QQ\,h^2+fg-h^2\bigr).
\leqno (3.8)
$$
Using (3.7) and (3.8), we can expand the RHS of (3.5) in the form
$\,\sum_{i=0}^4\,\QQ^i \,a_i$.  By quick inspection, $\,a_4=0$, $\,a_0=k^2$, 
and $\,a_1=sk-s^{\,\prime}k=(s-s^{\,\prime})\,k=-k^2$.  Also, by using
the definition of $\,s\,$ and $\,s^{\,\prime}$, we compute easily that
$$
a_2=(bc-d^2)\,(fg-h^2)-s^{\,\prime}s=(bh-df)\,(dg-ch)=r^{\,\prime}r,
\;\;\,\mbox{and} \leqno (3.9)
$$
\vspace{-10 mm}
$$
a_3=d^2(fg-h^2)+h^2(bc-d^2)-dh\,(s+s^{\,\prime})=(df-bh)\,(dg-ch)
=-r^{\,\prime}r. \leqno (3.10)
$$
With these computations of the $\,a_i$'s, the RHS of (3.5) becomes
\begin{eqnarray*}
- \QQ^3\, r^{\,\prime}r+\QQ^2\,r^{\,\prime}r-\QQ\,k^2+k^2
 \!&=&\! -\QQ^2\,r^{\,\prime}r\,\bigl( \QQ-1\bigr)-k^2\bigl( \QQ-1\bigr)\\
 \!&=&\! -q\cdot \bigl( \QQ^2\,r^{\,\prime}r+k^2\bigr),
\end{eqnarray*}
which is precisely the LHS of (3.5) as computed in (3.6).\qed

\bs 
Note that, while Thm.~3.4 generalizes Prop.~3.1, the proof of
the former is independent of that of the latter. Thus, in a mathematical
sense, we could have completely dispensed with Prop.~3.1. However, this
Proposition has clearly played an important role in the discovery (and 
formulation) of the formula (3.5), so we have included it for motivational
reasons.  In this regard, one might hope that there is also a shorter 
(or quicker) proof for (3.5) based on using a similar reduction (as in 
Remark 3.3) to the case where $\,X\,$ is {\it diagonal\/}. Unfortunately, 
the $\,q$-trace of a matrix is {\it not\/} invariant under (ordinary) 
matrix conjugation, so the standard results in the invariant theory 
of matrices (as developed in [Pr$_1$, Pr$_2$]; see also [Fo] and [KP]) 
does not apply directly to the setting of Thm.~3.4.  It is conceivable 
that some suitable form of a ``quantum invariant theory'' (based on an 
appropriate notion of quantum conjugation) might enable us to make the 
above reduction.  However, as far as we know, such a quantum invariant 
theory of matrices is not yet available.

\mk

\bs\nt
{\bf \S4. \ Quantum-Trace Determinantal Formula}

\bs
After our preliminary investigations on the case of $\,q$-traceless
matrices in \S3, we are now in a good position to give the full 
statement for the {\it quantum-trace determinantal formula.}  With 
respect to a fixed element $\,q\in S$, this formula expresses 
$\,q\cdot {\rm det}\,[X,Y]\,$ in terms of various quantities (including 
the $\,q$-traces) associated with the matrices $\,X,\,Y\,$ and $\,XY,\,YX$, 
grouped in descending powers of $\,\QQ:=1+q\in S$.  If $\,q\in S\,$
happens to be a unit (e.g.~a root of unity), we can then invert $\,q\,$ 
and get an equation just for $\,{\rm det}\,[X,Y]$.  The full statement 
is as follows.

\bs\nt
{\bf Theorem 4.1.} {\it Let $\,X,\,Y\in R={\mathbb M}_2(S)$, with
determinants $\,\dd,\,\dd^{\,\prime}$, traces $\,t,\,t^{\,\prime}$,
and $\,q$-traces $\,\tau,\,\tau^{\,\prime}$.  Also, let $\,\sigma
={\rm tr}_q(XY)$, and $\,\sigma^{\,\prime}={\rm tr}_q(YX)$.  Then}
$$
\begin{array}{c}
q\cdot {\rm det}\,[X,Y]=\QQ^2\,\dd^{\,\prime}\,\dd
-\QQ\,\bigl(\dd\,t^{\,\prime}\tau^{\,\prime}+\dd^{\,\prime}t\,\tau\bigr)
\quad \quad \quad \quad \quad \\
\quad \quad \quad \quad \quad \quad  \quad \quad \quad
+ \bigl( \dd\,\tau^{\,\prime\,2}+\dd^{\,\prime}\tau^2
+ {\rm tr}\,(XY)\,\tau^{\,\prime}\tau - \sigma^{\,\prime}\sigma\bigr).
\end{array}
\leqno (4.2)
$$
{\it Here, in the last parenthetical expression on the\/} RHS, 
{\it the first three terms constitute a quadratic form in 
$\,\tau^{\,\prime}\,$ and $\,\tau$, with coefficients 
$\,\dd,\,\dd^{\,\prime}$, and $\,{\rm tr}\,(XY)$.}

\bs
Before we proceed, let us first give the appropriate interpretations of 
the above formula in the two important special cases where $\,q={\pm 1}$. 
This is conceptually an important step, since Theorem 4.1 could hardly 
have come into existence without having had these two crucial special 
cases as its precursors.

\bs
{\it We first consider the case $\,q=1$.}  Here, $\,\tau=t$,
$\,\tau^{\,\prime}=t^{\,\prime}$, and $\,\sigma=\sigma^{\,\prime}
={\rm tr}\,(XY)$. (These are all ordinary traces.) In this case, after 
a simple combination of terms, (4.2) simplifies to the following 
\,{\bf trace version\/}\, of the determinantal formula\,:
$$
{\rm det}\;[X,Y] = 4\,\dd^{\,\prime}\dd - (\,{\rm tr}\,(XY)\,)^2
- \dd\,t^{\,\prime\,2} - \dd^{\,\prime}\,t^2 
+ {\rm tr}\,(XY)\;t^{\,\prime}t\,.
\leqno (4.3)
$$

\mk\nt
{\bf Remark 4.4.} (A) As in the general case, the last three terms in the 
formula above constitute a quadratic form in $\,t^{\,\prime}\,$ and $\,t$, 
with coefficients $\,-\dd,\,-\dd^{\,\prime}$, and $\,{\rm tr}\,(XY)$.  In 
the traceless case, this quadratic form drops out, and the formula boils
down to (3.2).  

\mk\nt
(B) Some more special cases of (4.3) are also worth noting. For instance, 
if $\,{\rm tr}\,(XY)=0$, (4.3) gives 
$\,{\rm det}\;[X,Y] = 4\,\dd^{\,\prime}\dd - \dd\,t^{\,\prime\,2} 
- \dd^{\,\prime}\,t^2$.  On the other hand, if $\,\dd=\dd^{\,\prime}=0$, 
(4.3) implies that $\,{\rm tr}\,(XY)\,$ divides $\,{\rm det}\,[X,Y]\,$ 
in the ring $\,S$.

\bs
From the viewpoint of Procesi's papers [Pr$_1$, Pr$_2$] (see also 
[KP:~\S2.4]), $\,{\rm det}\,[X,Y]\,$ should be expressible (in case 
$\,2\in {\rm U}(S)$) as a polynomial in the traces and determinants of 
$\,X,Y$, along with $\,{\rm tr}\,(XY)$. Thus, the existence of the 
formula (4.3) is entirely to be expected.  In the invariant theory of 
$\,2\times 2\,$ matrices, the two sides of the equation (4.3) represent 
the {\it Formanek element\/} associated with $\,X\,$ and $\,Y$; see, 
for instance, [JLS].  Note that, although the Formanek element in 
[JLS:~\S3.3] is expressed with the integer $\,4\,$ in the denominator, 
this ``$4$'' is eventually ``cancelled out'' to give the RHS of (4.3).  
We can, of course, also get a formula for $\,{\rm det}\,[X,Y]\,$ 
{\it using trace elements alone,} by expressing all determinants in 
(4.3) in terms of traces via the formula 
$\,2\;{\rm det}\,(M)=({\rm tr}\,(M))^2-{\rm tr}\,(M^2)\,$ (for every
$\,M\in {\mathbb M}_2(S)$).  This determinant elimination process 
results in the following equation:
$$
{\rm det}\;[X,Y]={\rm tr}\,(X^2)\,{\rm tr}\,(Y^2)
- \bigl({\rm tr}\,(XY)\bigr)^2 
- \frac{{\rm tr}\,(X^2)\,t^{\,\prime\,2} + {\rm tr}\,(Y^2)\,t^2}{2}
+ {\rm tr}\,(XY)\,t^{\,\prime}t\,.
$$
This time, the denominator ``$2$'' is no longer avoidable, so this
formula would be meaningful only over the rings $\,S\,$ in which 
$\,2\,$ is invertible. It is, in retrospect, rather fortunate that
the formula (4.3) is applicable to all commutative rings $\,S$.

\bs
{\it Next, we consider the case $\,q=-1\,$ in\/} Theorem 4.1. Here, 
for any $\,M=(m_{ij})\in {\mathbb M}_2(S)$, $\,{\rm tr}_q(M)\,$ is the 
{\it supertrace\/} $\,{\rm str}\,(M):=m_{11}-m_{22}\,$ defined in the
Introduction.  Thus, $\,\tau={\rm str}\,(X)$, $\,\tau^{\,\prime}
={\rm str}\,(Y)$, and $\,\sigma={\rm str}\,(XY)$, $\,\sigma^{\,\prime}
={\rm str}\,(YX)$.  Now for $\,q=-1$, we have $\,\QQ=1+(-1)=0$, so
all positive powers of $\,\QQ\,$ can be dropped\,! This leads to the
following remarkably simple \,{\bf supertrace version\/} of (4.2):
$$
-\,{\rm det}\;[X,Y] = \dd\,\tau^{\prime\,2} + \dd^{\,\prime}\,\tau^2
+{\rm tr}\,(XY)\;\tau^{\,\prime}\tau-{\rm str}\,(XY)\;{\rm str}\,(YX).
\leqno (4.5)
$$

\sk\nt
Here, as in Thm.~4.1, the first three terms constitute a quadratic form 
in $\,\tau^{\,\prime}\,$ and $\,\tau$, with coefficients 
$\,\dd,\, \dd^{\,\prime}$, and $\,{\rm tr}\,(XY)$. In (4.5), we have
chosen to keep the $\,(-1)$-factor on the LHS, to remind ourselves of 
the fact that this LHS is really $\,q\cdot {\rm det}\,[X,Y]$.

\bs\nt
{\bf Remark 4.6.} There are several ways to extend the definition of 
``str'' to higher matrix algebras. For instance, one may define ``str''
on $\,T={\mathbb M}_{2n}(S)\,$ by thinking of any $\,X\in T\,$ as a 
$\,2\times 2\,$ block matrix $\,(x_{ij})\,$ with four $\,n\times n\,$ 
blocks, and taking $\,{\rm str}\,(X)\,$ to be 
$\,{\rm tr}\,(x_{11})-{\rm tr}\,(x_{22})\in S$.  With this particular
definition of ``str'', the formula (4.5) would remain meaningful.  
Unfortunately, {\it it will no longer be true if $\,n>1$.}  In fact, let 
$\,X=Y={\rm diag}\,(I_n,0_n)$.  Then $\,[X,Y]=0$, so the LHS of (4.5) 
is zero.  Now $\,{\rm det}\,(X)={\rm det}\,(Y)=0$, $\,{\rm tr}\,(XY)=n$, 
and each of $\,X,\,Y,\,XY,\,YX\,$ has also supertrace $\,n$.  Thus, 
the RHS of (4.5) is $\,n^3-n^2=n^2(n-1)$.  So in this example, (4.5) 
{\it holds if and only if} $\;n=1\,$ (assuming, say, 
$\,{\rm char}\,(S)=0$). Exactly the same remark could have been made 
about the trace version of the determinantal formula in (4.3).

\bs
The following are some easy consequences of the supertrace formula (4.5).

\bs\nt
{\bf Corollary 4.7.} {\it Keep the notations used in the formula $(4.5)$.}

\mk
(1) {\it If $\,XY\,$ has a zero diagonal, then}
$\,-{\rm det}\;[X,Y]=\dd\,\tau^{\,\prime\,2}+\dd^{\,\prime}\,\tau^2$.

\sk
(2) {\it If $\,X\,$ has a constant diagonal, then} $\,-{\rm det}\;[X,Y]
=\dd\,\tau^{\,\prime\,2}-{\rm str}\,(XY)\,{\rm str}\,(YX)$.

\sk
(3) {\it If $\;XY=YX$, then} $\,\bigl({\rm str}\,(XY)\bigr)^2
=\dd\,\tau^{\,\prime\,2} + \dd^{\,\prime}\,\tau^2
+{\rm tr}\,(XY)\,\tau^{\,\prime}\tau$.

\bs\nt
{\bf Corollary 4.8.} {\it For any $\,a,b,c,d\in S$, we have}
$$
(ac-bd)^2 = ab\,(c-d)^2+cd\,(a-b)^2+(ac+bd)\,(a-b)\,(c-d).\leqno (4.9)
$$
{\it Alternatively, if we write $\,\tau=a-b\,$ and $\,\tau^{\,\prime}
=c-d$, then}
$$
(c\,\tau + b\,\tau^{\,\prime})^2=ab\,\tau^{\,\prime\,2}+cd\,\tau^2
+ (ac+bd)\,\tau^{\,\prime}\tau. \leqno (4.10)
$$

\sk\nt
{\bf Proof.}  The universal quartic identity (4.9) is just the result 
(4.7)(3), applied to a pair of {\it diagonal\/} matrices 
$\,X={\rm diag}\,(a,b)\,$ and $\,Y={\rm diag}\,(c,d)$.  The alternative 
form (4.10) follows from (4.9) upon noting that 
$\,c\,\tau+b\,\tau^{\,\prime}=c\,(a-b)+b\,(c-d)=ac-bd$. \qed

\bs
The quaternary quartic identity (4.9) does not seem well known. A search 
of the literature and standard websites such as [Pi] did not turn up 
this curious algebraic identity.  Of course, there is also a more 
sophisticated ``$\,q$-version'' of (4.9) (in five variables, including
$\,q\,)$, obtained by writing down the quantum-trace determinantal formula 
(4.2) for the diagonal matrices $\,X\,$ and $\,Y\,$ in the proof above.
Unlike (4.9), however, this quinary identity is no longer homogeneous.

\bs\nt
{\bf Corollary 4.11.} {\it Assume $\,S\,$ is a field, and keep the 
notations used in the formula $(4.5)$.  Then $\,[X,Y]\notin 
{\rm GL}_2(S)\,$ iff} $\,\dd\,\tau^{\prime\,2} + \dd^{\,\prime}\,\tau^2
+{\rm tr}\,(XY)\;\tau^{\,\prime}\tau = {\rm str}\,(XY)\;{\rm str}\,(YX)$.
{\it Also, $\,X\,$ is a scalar matrix iff this equation holds for all 
$\,Y\in {\mathbb M}_2(S)$.}

\bs\nt
{\bf Proof.} The first statement follows from (4.5), since 
$\,{\rm det}\;[X,Y]=0\,$ iff $\,[X,Y]\notin {\rm GL}_2(S)$.  For the
second statement, the ``only if'' part is clear, and the ``if'' part
follows from the easy fact (see [GL], or [KL$_2$:~(5.14)]) that, for 
any non-scalar matrix $\,X\,$ (over a field $\,S$), there exists some 
$\,Y\in {\mathbb M}_2(S)\,$ such that $\,[X,Y]\in {\rm GL}_2(S)$. \qed

\bs
We shall now begin to work toward the proof of Thm.~4.1.  Of course, 
once the determinantal formula (4.2) is written down {\it explicitly,} 
a direct check on {\tt Singular} or {\tt Macaulay2} will instantly 
confirm that it is a universal identity for polynomials in nine 
commuting variables (the symbol $\,q\,$ together with the eight entries 
of $\,X\,$ and $\,Y$).  However, such a machine checking exercise would 
reveal no reason whatsoever for the truth of the formula.  In view of 
this, we feel it still imperative to give a detailed conventional 
mathematical proof for the formula (4.2).  Our proof is preceded by
the following lemma.

\bs\nt
{\bf Lemma 4.12.} {\it For $\,X,Y\in {\mathbb M}_2(S)\,$ with notations 
as in\/} Theorem 4.1, {\it we have}
$$
\sigma + \sigma^{\,\prime}-t^{\,\prime}\tau-t\,\tau^{\,\prime}
= \QQ\,\bigl(\,{\rm tr}\,(XY)-t^{\,\prime}t\,\bigr). \leqno (4.13)
$$

\sk\nt
{\bf Proof.} By working generically, we may assume, as in the proof of
Prop.~2.2, that $\,2\,$ is invertible in $\,S$. In this case, we can
write $\,{\rm det}\;X=[\,({\rm tr}\;X)^2-{\rm tr}\;(X^2)]/2$.
Substituting this into the Cayley-Hamilton equation
$\,X^2-({\rm tr}\;X)\,X+({\rm det}\;X)\,I_2=0$, and polarizing the
resulting equation (via $\,X\mapsto X+Y$), we get
$$
XY+YX-({\rm tr}\;Y)\,X-({\rm tr}\;X)\,Y 
+ ({\rm tr}\;X)\,({\rm tr}\;Y)\,I_2 - {\rm tr}\;(XY)\,I_2=0.
$$
Taking $\,q$-traces on both sides gives the desired equation (4.13).\qed

\bs
Before proceeding to the proof of Theorem 4.1, we record a couple of
consequences of the Lemma above.  

\bs\nt
{\bf Corollary 4.14.} (1) {\it For any traceless matrices $\,X,Y\in 
{\mathbb M}_2(S)$, we have}
$$
{\rm tr}_q(XY)+{\rm tr}_q(YX) = \QQ\;{\rm tr}\,(XY). \leqno (4.15)
$$
\sk\nt
(2) {\it For any matrices $\,X,Y\in {\mathbb M}_2(S)$, we have}
$$
{\rm str}\,(XY)+{\rm str}\,(YX)={\rm tr}\,(X)\,{\rm str}\,(Y)
+{\rm tr}\,(Y)\,{\rm str}\,(X).\leqno (4.16)
$$
{\it In particular, $\,{\rm str}\,(YX)=-{\rm str}\,(XY)\,$ if $\,X$
has a zero diagonal\/, or $\,X,Y$ are both traceless, or $\,X,Y\,$ are
both supertraceless.}

\bs\nt
{\bf Proof.} (1) is obtained from Lemma 4.12 by setting 
$\,t=t^{\,\prime}=0$. [\,In this special case, (4.15) expresses the fact 
that the ordinary trace $\,{\rm tr}\,(XY)\,$ is a ``quantum average'' 
of the quantum traces $\,{\rm tr}_q(XY)\,$ and $\,{\rm tr}_q(YX)\,$ 
(for any given $\,q\in S$).]

\mk\nt
(2) is obtained by specializing Lemma 4.12 to the case $\,q=-1$, where
the $\,q$-trace becomes the supertrace.  [\,In addition, we can check
easily that both sides of (4.16) are equal to $\,2\;{\rm str}\,(X\ast Y)$, 
where $\,X\ast Y\,$ is the Hadamard product of $\,X\,$ and $\,Y\,$ 
(obtained by ``entry-wise multiplication'' of the two matrices).]\qed

\bs
We have now all the necessary tools with which to verify our formula (4.2).

\bs\nt
{\bf Proof of Theorem 4.1.} To begin with, note that if the two matrices
$\,X,\,Y\,$ are both $\,q$-traceless (that is, $\,\tau=\tau^{\,\prime}=0$),
then all terms on the RHS of the determinantal formula (4.2) drop out ---
except the first and the last terms.  In this case then, we know that 
(4.2) holds, thanks to Thm.~3.4.  To prove (4.2) in general, {\it we
should then try to make a reduction to the $\,q$-traceless case.}
As before, we may assume that $\,S\,$ is the free commutative 
$\,{\mathbb Q}$-algebra generated by $\,q\,$ and the eight entries of 
$\,X\,$ and $\,Y$. In particular, $\,\QQ^{-1}\,$ exists in the quotient
field of the integral domain $\,S$.  To make the desired reduction, let
$\,X_0=X-\QQ^{-1}\tau\,I_2\,$ and $\,Y_0=Y-\QQ^{-1}\tau^{\,\prime} I_2$. 
These have $\,q$-traces zero since $\,\tau,\,\tau^{\,\prime}\,$ were the
$\,q$-traces of $\,X\,$ and $\,Y$, and obviously $\,[X,Y]=[X_0,Y_0]$.  
Therefore, by Thm.~3.4, we have
$$
q\cdot {\rm det}\,[X,Y]=\QQ^2\,{\rm det}\,(X_0Y_0)
-{\rm tr}_q(X_0Y_0)\,{\rm tr}_q(Y_0X_0).\leqno (4.17)
$$
Our job now is to compute the RHS of (4.17) in terms of the various
quantities associated with $\,X,Y,XY$, and $\,YX$. Taking $\,q$-traces
on the equation
$$
X_0Y_0=XY-\QQ^{-1}(\tau^{\,\prime}X+\tau\,Y)
+\QQ^{-2}\tau^{\,\prime}\tau\,I_2\,, \leqno (4.18)
$$
we see that $\,{\rm tr}_q(X_0Y_0)=\sigma-\QQ^{-1}\tau^{\,\prime}\tau\,$ 
where $\,\sigma={\rm tr}_q(XY)$, and similarly,
$\,{\rm tr}_q(Y_0X_0)=\sigma^{\,\prime}-\QQ^{-1}\tau^{\,\prime}\tau$,
where $\,\sigma^{\,\prime}={\rm tr}_q(YX)$.  On the other hand, since 
$\,{\rm det}\,\bigl(X-\lambda\,I_2\bigr)=\lambda^2-t\,\lambda+\dd\,$ for 
any parameter $\,\lambda\in S\,$ (notations as in Thm.~4.1), we have,
for $\,\lambda=\QQ^{-1}\tau$\,:
$$
{\rm det}\,(X_0)=\QQ^{-2}\tau^2-\QQ^{-1}t\,\tau+\dd\,; \;\mbox{similarly,}
\;\;{\rm det}\,(Y_0)=\QQ^{-2}\tau^{\,\prime\,2}-\QQ^{-1}t^{\,\prime}
\tau^{\,\prime}+\dd^{\,\prime}. 
$$
Substituting all of these expressions into the RHS of (4.17), we get
$$
\begin{array}{c}
q\cdot {\rm det}\,[X,Y]=
\bigl(\,\QQ\,\dd-t\,\tau + \QQ^{-1}\tau^2 \bigr)
\bigl(\,\QQ\,\dd^{\,\prime}-t^{\,\prime}\tau^{\,\prime}
+ \QQ^{-1}\tau^{\,\prime\,2} \bigr) \\
\;\;\;\;\;\;\;\;\;\;\;\;\;\;\;\;\;\,
-\bigl( \sigma-\QQ^{-1}\tau^{\,\prime}\tau\bigr)
\bigl( \sigma^{\,\prime}-\QQ^{-1}\tau^{\,\prime}\tau \bigr).
\end{array}  \leqno (4.19)
$$
Expanding the RHS {\it formally\/} into $\,\sum_{i=-2}^2\,\QQ^i\,b_i$, 
we have clearly $\,b_{-2}=0$, $\,b_2=\dd^{\,\prime}\dd$, and 
$\,b_1=-\bigl(\dd\,t^{\,\prime}\tau^{\,\prime}+\dd^{\,\prime}t\,\tau\bigr)$.
Thus, $\,\QQ^{2}\,b_2+\QQ\,b_1\,$ already produces the first two groups of 
terms on the RHS of the formula (4.2).  The remaining terms on the RHS 
of (4.19) are
\begin{eqnarray*}
b_0+\QQ^{-1}b_{-1} \!&=&\!
\bigl( \dd\,\tau^{\,\prime\,2}+\dd^{\,\prime}\tau^2
+ t^{\,\prime}t\,\tau^{\,\prime}\tau - \sigma^{\,\prime}\sigma \bigr) 
+ \QQ^{-1} \tau^{\,\prime}\tau\,\bigl(\sigma + \sigma^{\,\prime}
- t^{\,\prime}\tau - t\,\tau^{\,\prime}\bigr) \\
\!&=&\!  \bigl( \dd\,\tau^{\,\prime\,2}+\dd^{\,\prime}\tau^2
+ t^{\,\prime}t\,\tau^{\,\prime}\tau - \sigma^{\,\prime}\sigma \bigr) 
+ \tau^{\,\prime}\tau \,\bigl(\,{\rm tr}\,(XY)-t^{\,\prime}t\,\bigr)
\end{eqnarray*}
in view of Lemma 4.12.  After cancelling the two 
$\,t^{\,\prime} t\, \tau^{\,\prime}\tau\,$ terms, we get precisely the 
last group of terms in the desired determinantal formula (4.2)\,!\qed

\bs\nt
{\bf Remark 4.20.} Of course, proving the quantum-trace version of the
determinantal formula (4.2) in one stroke for all $\,q\,$ makes it 
unnecessary, for instance, to handle separately the cases $\,q=1\,$ and
$\,q=-1$. But more discerningly, working directly in the quantum-trace 
case actually makes the proof of (4.2) {\it easier\/} as it enables us 
to ``manage'' many terms at once by organizing (and simplifying) them in 
``$[\,2\,]$-adic expansions'' (as in the proof above).  The same proof, 
written out in the special cases $\,q=1\,$ or $\,q=-1\,$ would look 
harder and more confusing since the pattern of the $[\,2\,]$-adic 
expansions would no longer be apparent.  The same remark could have  
been made about the proof of Thm.~3.4.

\mk

\bs\nt
{\bf \S5. \ Relations to Binary Quadratic Forms}

\bs
The last two sections of this paper are devoted to some applications 
of the two determinantal formulas obtained in (4.3) and (4.5). The first
applications, given in this section, offer a characteristic-free
generalization of a theorem of Olga Taussky relating the determinants of 
$2\times 2$ commutators of integral matrices to norms in quadratic 
extensions of $\,{\mathbb Q}\,$, and some extensions of this theorem to 
the setting of matrices over commutative rings.

\bs
In [Ta$_1$], Taussky showed that, {\it if $\,X,Y\in 
{\mathbb M}_2({\mathbb Z})\,$ and an eigenvalue $\,\omega\,$ of $\,X\,$ 
is irrational\/}, {\it then $\,-{\rm det}\;[X,Y]\,$ is a norm from the 
quadratic number field $\,{\mathbb Q}\,(\omega)$.} A converse of this 
theorem was obtained in [Ta$_2$], where Taussky proved that, 
{\it if $\,n\in {\mathbb Q}\,$ is a norm from a quadratic number 
field $\,K$, then $\,n=-{\rm det}\;[X,Y]\,$ for some $\,X,Y\in 
{\mathbb M}_2({\mathbb Q})\,$ such that $\,X\,$ has its eigenvalues 
in $\,K$.}  Although Taussky assumed that $\,X,Y\,$ were integral 
matrices in the first theorem above, this assumption was not needed, 
so both of her theorems may be thought of as results on rational 
matrices.  Actually, the use of the rational field $\,{\mathbb Q}\,$ 
is also not crucial, so one may try to further replace $\,{\mathbb Q}\,$ 
by a field $\,F$.  However, Taussky's proofs in [Ta$_1$, Ta$_2$]
(and even her later proof using cyclic algebras in [Ta$_3$]) assumed 
implicitly that $\,{\rm char}\,(F)\neq 2$, and did not apply
to all fields.

\bs
In the first half of this section, we shall present a new view of
both of Taussky's results, formulating them as a ``commutator 
characterization'' for the norm elements under any quadratic field 
extension $\,K/F\,$ (separable or otherwise).  Here, we are able to 
give a rather short proof (motivated by the determinantal formula 
(4.3)) {\it that works uniformly in all characteristics,\/} and is 
completely within the realm of matrix theory (independently of the
splitting criterion for cyclic algebras used in [Ta$_3$]).  Furthermore, 
the proofs of the ``if'' part and the ``only if'' part below are based 
essentially on one single argument, and the ``if'' part will, later 
in the section, lead to a constructive generic version of the same 
result for commutative rings.

\bs\nt
{\bf Taussky's Norm Theorem 5.1.} {\it Let $\,K/F\,$ be any quadratic 
extension of fields of any characteristic.  Then an element $\,n\in F\,$ 
is a norm from $\,K\,$ iff $\,n=-{\rm det}\;[X,Y]\,$ for some
$\,X,Y\in {\mathbb M}_2(F)\,$ such that $\,K\,$ is the splitting field 
of the characteristic polynomial of $\,X\,$.}

\bs\nt
{\bf Proof.} We'll first prove the harder ``if'' part. Given 
$\,n=-{\rm det}\;[X,Y]\,$ as in the theorem, let $\,t={\rm tr}\,(X)$,  
$\,\delta={\rm det}\,(X)$, and let $\,\omega\,$ be an eigenvalue of 
$\,X$.  By assumption, $\,K=F(\omega)$, so the minimal polynomial for
$\,\omega\,$ over $\,F\,$ is $\,f(\lambda)=\lambda^2-t\,\lambda+\delta\,$
(the characteristic polynomial of $\,X$).  With respect to the 
$\,F$-basis $\,\{1,\omega\}\,$ on $\,K$, the norm form of $\,K/F\,$ is 
easily computed to be\footnote{If we had used $\,\{1,-\omega\}\,$ as 
basis instead, the norm form would have been $\,x^2-t\,xy+\delta\,y^2$,
which is precisely the homogenization of the characteristic polynomial 
of $\,X$.}
$$
{\bf N}\,(x,y):={\rm N}_{K/F}(x+y\,\omega)=x^2+y\,(t\,x+\delta\,y),\;\;\,
\mbox{where}\;\;x,y\in F.\leqno (5.2)
$$
Since $\,f(\lambda)\,$ is irreducible over $\,F$, we may assume (after a 
conjugation in $\,{\mathbb M}_2(F)$) that $\,X\,$ is in its rational  
canonical form; that is, $\,X=\begin{pmatrix} 0&-\dd \\1&t\end{pmatrix}$.  
Also, after subtracting a scalar matrix from $\,Y\,$ (which does not 
change $\,[X,Y]\,$), we may assume that $\,Y=\begin{pmatrix} 
e&f\\g&0\end{pmatrix}\,$ (for some $\,e,f,g\in F$). 
Since $\,{\rm tr}\,(XY)=f-\dd\,g$, the trace formula (4.3) yields:
\begin{eqnarray*}
-{\rm det}\;[X,Y] &=& 4\,\dd fg+(f-\dd\,g)^2
+ \dd\,e^2 - fg\,t^2 - (f-\dd\,g)\,t\,e   \\
&=& (f+\dd\,g)^2+(e+t\,g)\,(\dd\,e-tf).
\end{eqnarray*}
Letting $\,x:=-(f+\dd\,g)\,$ and $\,y:=e+t\,g$, we get
$\,-{\rm det}\;[X,Y]=x^2+y\,(\delta\,e-tf)$.  (Of course, 
this negative determinant could also have been gotten from the supertrace
formula (4.5), or even from a direct determinant computation.\footnote{In 
fact, if we use the supertrace formula (4.5), $\,x=-(f+\dd\,g)\,$ will 
show up naturally as $\,{\rm str}\,(XY)$, and we will have
$\,{\rm str}\,(YX)=f+\dd\,g=-x$, as is also predicted by the formula 
(4.16).})  Noting that
$$
t\,x+\dd\,y=-t\,(f+\dd\,g)+\dd\,(e+t\,g)
=\dd\,e-tf,\leqno (5.3)
$$ 
we conclude from (5.2) that $\,n=-{\rm det}\;[X,Y]={\bf N}\,(x,y)\in
{\rm N}_{K/F}(K)$.

\mk
The converse is now easy! Indeed, we can completely bypass the work in 
[Ta$_2$], and simply ``reverse'' the above argument to get what we want,
as follows. Let $\,n\in F\,$ be a norm from $\,K$. Write $\,K=F(\omega)\,$ 
for a primitive element $\,\omega$, and let $\,\lambda^2-t\,\lambda+\delta\,$ 
be the minimal polynomial of $\,\omega\,$ over $\,F$. Then 
$\,n={\rm N}_{K/F}(x+y\,\omega)\,$ for some $\,x,y\in F$. Defining 
$\,X:=\begin{pmatrix} 0&-\dd \\1&t\end{pmatrix}\,$ and $\,Y=\begin{pmatrix} 
y&-x\\0&0\end{pmatrix}$, the computation in the last paragraph (with 
$\,g=0,\,f=-x$, and $\,e=y$) gives 
$\,-{\rm det}\;[X,Y]={\rm N}_{K/F}(x+y\,\omega)=n$.  Of course, 
the splitting field of the characteristic polynomial of $\,X\,$ is just
$\,K$. \qed

\bs\nt
{\bf Remark 5.4.} The following observation on the ``if'' part of
Theorem 5.1 is in order.  Let $\,F={\mathbb Q}\,$, and assume (as in
[Ta$_1$]) that $\,X,Y\in {\mathbb M}_2({\mathbb Z})$. If $\,X\,$ is in 
its rational canonical form $\,\begin{pmatrix} 0&-\dd \\1&t\end{pmatrix}$, 
the proof of the ``if'' part above shows that $\,-{\rm det}\;[X,Y]\,$ is, 
in fact, the norm of an {\it algebraic integer\/} in the quadratic field 
$\,K$.  However, in general, this need not be the case, as was pointed 
out by Taussky in [Ta$_2$:~p.\,1]. For an explicit example, take 
$\,X=\begin{pmatrix}0&4\\-2&1\end{pmatrix}\,$ and $\,Y:=\begin{pmatrix}
4&3\\3&0\end{pmatrix}$.  Here, the quadratic field $\,K\,$ in question 
is $\,{\mathbb Q}\,(\sqrt{-31}\,)$, and $\,-{\rm det}\;[X,Y]=419\,$ is 
not the norm of an algebraic integer from $\,K$, since $\,x^2+31\,y^2
=2^2\cdot 419=1676\,$ has no solution in $\,{\mathbb Z}$.  (Of course, 
$\,X\,$ is only ``close'' --- but not equal --- to its rational canonical 
form!)  Nevertheless, in confirmation of the ``if'' part of Thm.~5.1, 
$\,4^2\cdot 419=6704=x^2+31\,y^2\,$ is solved by 
$\,(x,y)=({\pm 77},{\pm 5})$, so $\,419={\rm N}_{K/{\mathbb Q}}(\alpha)\,$ 
for $\,\alpha=\bigl({\pm 77}\,{\pm}\,5\sqrt{-31}\,\bigr)/4\,$ in 
$\,{\mathbb Q}\,(\sqrt{-31}\,)$.  For more information on this example, 
see Remark 5.15(B) below.

\bs
Since (5.1) was formulated as a field-theoretic theorem, a natural 
question to ask would be whether something in a similar spirit can be 
said about commutative rings.  Note that the proof for the ``if'' part 
of Theorem 5.1 does not extend to rings, since we can no longer apply 
the standard linear algebra theorem on rational canonical forms.  Of 
course, the field-theoretic theorem, applied to the quotient field 
of a suitable generic ring, would give an existential norm formula 
on $\,-{\rm det}\;[X,Y]\,$ for a pair of generic matrices $\,X,\,Y$. 
However, this formula would involve an unknown ``denominator'' factor, 
which would make it only a strictly formal result.  A useful
``ring-theoretic version'' of Theorem 5.1 should thus be one that gives 
an {\it implementation\/} of such a formula, with explicit information 
on the denominator factor.  Such a version will be given (for 
any commutative ring) in Theorem 5.10 below, where we'll show that 
{\it the denominator factor can actually be taken to be either one of 
the off-diagonal entries of the matrix $\,X$.}  The proof of this 
theorem is based on a further exploitation of the explicit determinant 
computation of commutators in the proof of the ``if'' part of 
Theorem 5.1.

\bs
As in the earlier sections, {\it $\,S\,$ will continue to denote a
commutative ring.}  Instead of working with the norms from various
degree $2$ extensions of $\,S$, we now choose to work directly 
with binary quadratic forms over $\,S$.  For any $\,s,t,\dd\in S$, 
let us denote the ``value set'' (over $\,S$) of the quadratic form
$\,s\,x^2+t\,xy+\dd\,y^2\,$  by $\,V[\,s,t,\dd\,]$; that is,
$$
V[\,s,t,\dd\,]:=\{\,s\,r_1^2+t\,r_1r_2+\dd\,r_2^2:\;r_1,r_2\in S\,\}.
\leqno (5.5)
$$
In the case where $\,t=0$, we'll simply write $\,V[\,s,\dd\,]\,$ for
$\,V[\,s,0,\dd\,]$.  Over the ring of integers, of course, the study 
of these value sets is an important and time-honored enterprise that 
goes back to the classical work of Fermat, Euler, Lagrange, Legendre, 
and Gauss. In the rest of this section, we'll work over a commutative 
ring $\,S\,$ in the case $\,s=1$, since $\,x^2+t\,xy+\dd\,y^2\,$ 
arises precisely as a norm form of the quadratic $\,S$-algebra 
$\,S\,[\,\lambda\,]/(\lambda^2-t\,\lambda+\dd)$. Before coming to the 
ring-theoretic version of Thm.~5.1, we first recall the following 
elementary result on the value sets $\,V[\,1,t,\dd\,]\,$ over $\,S$.
A short proof is included for the reader's convenience.

\bs\nt
{\bf Proposition 5.6.} {\it For any $\,t,\,\dd\in S\,$ and 
$\,\Delta:=t^2-4\,\dd$, we have the following inclusions\/}:
$$
4\,V[\,1,t,\dd\,]\subseteq V[\,1,-\Delta\,]
\subseteq V[\,1,t,\dd\,].  \leqno (5.7)
$$
{\it If $\;2\,$ is invertible in $\,S$, then 
$\,V[\,1,-\Delta\,]=V[\,1,t,\dd\,]$.}

\bs\nt
{\bf Proof.}  For any $\,w,z\in S$, we have an identity:
$$
w^2-(t^2-4\,\dd)\,z^2=(w-t\,z)^2+t\,(w-tz)\,(2z)+\dd\,(2z)^2.\leqno (5.8)
$$
This implies that the set $\,\{x^2+txy+\dd\,y^2:x\in S,\,y\in 2S\}\,$
is equal to the set $\,\{w^2-\Delta\,z^2:w,z\in S\}$.  The second inclusion
in (5.7) follows from this observation.  The first inclusion follows from
the usual ``completion of squares'' identity:
$$
4\,(x^2+t\,xy+\dd\,y^2)=(2\,x+t\,y)^2-(t^2-4\,\dd)\,y^2.\leqno (5.9)
$$
The last conclusion of the Proposition is clear from (5.7).\qed

\bs
We are now in a position to extend Theorem 5.1 to the setting of rings.

\bs\nt
{\bf Norm Theorem 5.10. (Ring Version)} (1) {\it Given $\,t,\,\dd\in S$, 
any $\,n\in V[\,1,t,\dd\,]\,$ has the form $\,-{\rm det}\;[X,Y]\,$ for 
some $\,X,Y\in {\mathbb M}_2(S)\,$ such that $\,{\rm tr}\,(X)=t\,$ and 
$\,{\rm det}\,(X)=\dd$.}

\sk\nt
\nt (2) {\it For any $\,X=\begin{pmatrix}a&b\\c&d\end{pmatrix}\in 
{\mathbb M}_2(S)$, let $\,t={\rm tr}\,(X)$, $\,\dd={\rm det}\,(X)$, 
and $\,\Delta=t^2-4\,\dd$.  For any $\,Y\in {\mathbb M}_2(S)$, we have 
$\,-c^2\,{\rm det}\;[X,Y]\in V[\,1,t,\dd\,]$, and 
$\,-4\,c^2\,{\rm det}\;[X,Y]\in V[\,1,-\Delta\,]$.}

\mk\nt
(3) {\it Keep the notations in $(2)$ above.  If $\;t={\rm tr}\,(X)\,$
and $\,t^{\,\prime}={\rm tr}\,(Y)\,$ are both in $\,2\,S$, then 
$\;-c^2\,{\rm det}\;[X,Y]\in V[\,1,-\Delta\,]$.}

\bs\nt
{\bf Proof.} (1) follows from the last paragraph in the proof of 
Theorem 5.1, since the construction there works over any commutative 
ring $\,S$. 

\sk\nt
(2) After subtracting a scalar matrix from $\,Y$, we may assume that $\,Y
=\begin{pmatrix}e&f\\g&0\end{pmatrix}$.  We can work {\it generically\/} 
and thus assume that $\,S\,$ is the polynomial ring over $\,{\mathbb Z}\,$ 
generated by the seven (commuting) variables $\,a,b,c,d,e,f,g$. In this 
way, $\,c^{-1}\,$ makes sense in the quotient field $\,F\,$ of $\,S$. 
Let $\,X_1=\begin{pmatrix}0&c^{-1}b\\1&c^{-1}(d-a)\end{pmatrix}\in 
{\mathbb M}_2(F)$.  Applying formally the calculation in the proof of 
the ``if'' part of Thm.~5.1, we can write
$$
-\,{\rm det}\;[X_1,Y]=x^2+c^{-1}(d-a)\,xy-c^{-1}b\,y^2,\leqno (5.11)
$$
where $\,x=-(f-c^{-1}b\,g)$, and $\,y=e+g\,c^{-1}(d-a)$. From these, 
we have $\,cx,\,cy\in S$.  Defining $\,X_2=c\cdot X_1=\begin{pmatrix}
0&b\\c&d-a\end{pmatrix}$ and multiplying (5.11) by $\,c^4$, we see that
$$
-\,c^2\,{\rm det}\;[X_2,Y]=(c^2x)^2+(d-a)\,(c^2x)\,(c\,y)-b\,c\,(c\,y)^2.
\leqno (5.12)
$$
Letting $\,\alpha=c^2x\in S\,$ and $\,\beta=cy\in S$, the RHS of (5.12)
can be transformed as follows:
\begin{eqnarray*}
\alpha^2+(d-a)\,\alpha\,\beta-b\,c\,\beta^2 &=& 
\alpha^2+(t-2\,a)\,\alpha\,\beta+(\dd-ad)\,\beta^2 \\
&=& (\alpha-a\,\beta)^2+t\,\alpha\,\beta-a\,(a+d)\,\beta^2+\dd\,\beta^2\\
&=& (\alpha-a\,\beta)^2+t\,(\alpha-a\,\beta)\,\beta+\dd\,\beta^2
\in V[\,1,t,\dd\,].
\end{eqnarray*}
Since $\,[X_2,Y]=[X_2+a\,I_2,Y]=[X,Y]$, this proves the first conclusion 
in (2). The second conclusion follows from this and the first inclusion 
in (5.7).

\mk\nt
(3) Write $\,t=2\,s\,$ and $\,t^{\,\prime}=2\,s^{\,\prime}\,$ (for 
suitable $\,s,s^{\,\prime}\in S$). {\it Suppose the desired conclusion
is true for traceless matrices.}  Then it holds for $\,X_0:=X-s\,I_2\,$ 
and $\,Y_0:=Y-s^{\,\prime}I_2$; that is, $\,-c^2\,{\rm det}\;[X_0,Y_0]
\in V[\,1,\,4\,{\rm det}\,(X_0)\,]$.  (Note that $\,X_0\,$ has 
discriminant $\,-4\,{\rm det}\,(X_0)$, and the $\,(2,1)$-entry of 
$\,X_0\,$ remains to be $\,c$.)  To compute $\,{\rm det}\,(X_0)$, we 
use the fact that $\,{\rm det}\,(X-\lambda\,I_2)=\lambda^2-t\,\lambda
+\delta$.  For $\,\lambda=s$, this leads to
$$
4\;{\rm det}\,(X_0)=4\,(s^2-t\,s+\dd\,)=t^2-2\,t^2+4\,\dd
=-(t^2-4\,\dd\,)=-\Delta\,. \leqno (5.13)
$$
Therefore, $\,-c^2\,{\rm det}\;[X_0,Y_0]\in V[\,1,-\Delta\,]$.
Since $\,[X,Y]=[X_0,Y_0]$, this proves (3).  Starting afresh, {\it we 
may thus assume that\/} $\,X=\begin{pmatrix}a&b\\c&-a\end{pmatrix}$ 
and $\,Y=\begin{pmatrix}e&f\\g&-e\end{pmatrix}$.  In this case, 
(3) can be proved by checking the explicit equation 
$\,-c^2\,{\rm det}\;[X,Y]=P^2-\Delta\,Q^2$, where
$$
P=2\,a\,(a\,g-c\,e)+c\,(b\,g-cf),\;\;\,\mbox{and}\;\;\;Q=a\,g-c\,e. 
\leqno (5.14)
$$
Since in any case such an equation can be quickly checked by hand or by 
machine, we will not give its detailed derivation here.\qed

\bs\nt
{\bf Remark 5.15.} (A) Note that, in the case where $\,S\,$ is a field
$\,F$, the results in (1) and (2) above do retrieve the Norm Theorem 5.1. 
In fact, in part (2), if the characteristic polynomial of the matrix 
$\,X\,$ has a quadratic splitting field $\,K/F$, then the off-diagonal 
entries $\,b,c\,$ of $\,X\,$ cannot both be zero.  If $\,c\neq 0$, then 
the conclusion $\,-c^2\,{\rm det}\;[X,Y]\in V[\,1,t,\dd\,]\,$ amounts 
to $\,-{\rm det}\;[X,Y]\in V[\,1,t,\dd\,]$, since $\,c\,$ is invertible, 
and $\,V[\,1,t,\dd\,]\,$ is closed under multiplication by squares. If 
$\,b\neq 0$, a simple transposition argument gives the same conclusion.  
But of course, the proof of Thm.~5.10 would not have been possible 
if we had not first worked out the proof of the field-theoretic version 
Thm.~5.1.

\mk\nt
(B) Since the proof of (5.10)(2) is completely constructive, we can
very easily implement it and test its accuracy.  For instance, let us 
apply it to the two matrices $\,X,\,Y\,$ in Remark 5.4 over the ring 
$\,S={\mathbb Z}\,$.  Here, $\,Y\,$ already has the desired form in the
proof of (5.10)(2), with $\,e=4\,$ and $\,f=g=3$, while $\,(a,b,c,d)
=(0,4,-2,1)$, with $\,t=1$, $\,\dd=8$, $\,\Delta=-31$, and 
$\,-{\rm det}\;[X,Y]=419$.  We know that $\,419\notin V[\,1,1,8\,]\,$ 
(since $\,419\,$ is not the norm of an algebraic integer in 
$\,{\mathbb Q}\,(\sqrt{-31}\,)$), {\it so the $\,c^2\,$ factor cannot 
be dropped from the first conclusion of\/} (5.10)(2).  On the other 
hand, following the proof of (5.10)(2), we compute easily that 
$\,\alpha=-36$, and $\,\beta=-5$.  Since $\,a=0\,$ and $\,c=-2$, this 
proof predicts that $\,c^2\cdot 419=1676\in V[\,1,1,8\,]$, 
with the equation $\,1676=\alpha^2+\alpha\,\beta+8\,\beta^2\,$ solved 
by $(\alpha,\beta)=(-36,-5)$. In view of this and the equation (5.9), 
the last part of (5.10)(2) also predicts that $\,4\,c^2\cdot 419=6704\in 
V[\,1,31\,]$, with the equation $6704=\alpha_0^2+31\,\beta_0^2\,$ 
solved by $\,(\alpha_0,\beta_0)=(2\,\alpha+\beta,\beta)=(-77,-5)$, as 
we have already mentioned in Remark 5.4.  Since $\,c^2\cdot 419=1676
\notin V[\,1,31\,]$, this shows that {\it the factor of $\,4\,$ also 
cannot be dropped from the second conclusion of\/} (5.10)(2).  (We
leave it to the reader to check the same statement if we had started
instead with $\,-{\rm det}\,[\,Y,X]=419$.  Note that (5.10)(3) does 
not apply to either case since $\,{\rm tr}\,(X)=1\,$ and
$\,{\rm tr}\,(Y)=4\,$ are {\it not both even\/}\,!\,) 

\bs
We'll close this section with a supplement to Theorem 5.10 in 
the case where the matrix $\,X\,$ has a constant diagonal. In this 
case, we have good control on the values of $\,-{\rm det}\,[X,Y]\,$ 
without pre-multiplying them by the factor $\,c^2$. The following
result is not covered by Theorem 5.10, but its proof is straightforward 
in light of the supertrace determinantal formula (applied in its special 
form in (4.7)(2)).

\mk\nt
{\bf Proposition 5.16.} {\it Let $\,X=\begin{pmatrix}a&b\\c&a\end{pmatrix}
\in R={\mathbb M}_2(S)$, and assume that $\,bS+cS=rS\,$ for some $\,r\in S$.
Then $\,\{-{\rm det}\;[X,Y]: Y\in R\}=V[\,r^2,-bc\,]$.  If $\,b,\,c\,$ are
coprime in $\,S$, this set is equal to $\,V[\,1,-bc\,]$.}

\mk\nt
{\bf Proof.} We may assume that $\,a=0$, and that $\,Y=\begin{pmatrix}
w&x\\y&0\end{pmatrix}$. Then, by (4.7)(2) (or by a direct computation), 
we have $\,-\,{\rm det}\,[X,Y]=(by-cx)^2-bcw^2$. Since $\,w\,$ 
ranges over $\,S\,$ and $\,by-cx\,$ ranges over the principal ideal 
$\,r\,S$, these values comprise precisely the set $\,V[\,r^2,-bc\,]
\subseteq V[\,1,-bc\,]$. If $\,bS+cS=S$, we can take $\,r=1$, in which 
case the inclusion becomes an equality. \qed

\mk

\bs\nt
{\bf \S6. \ Applications to Matrix Factorizations and Affine Curves}

\bs
Continuing the work in \S5, we shall give in this section some applications 
of the {\it supertrace\/} determinantal formula (4.5).  The main themes 
of our study will now be the factorization of $2\times 2$ matrices, and 
the solution of certain quadratic diophantine equations over a commutative 
ring $\,S$. The norm forms of quadratic ring extensions over $\,S\,$
studied in the last section are the binary quadratic forms $\,x^2+t\,xy
+\delta\,y^2$, which are {\it monic\/} in $\,x$.  In this section, we 
shall take up the case of a binary {\it diagonal\/} quadratic form 
$\,p\,x^2+q\,y^2\,$ (which is no longer monic in $\,x$).  In the spirit 
of the results (5.1) and (5.10), we would like to give a ``commutator 
characterization'' for the value set of such a diagonal form over $\,S$; 
that is, 
$$
V[\,p,q\,]:=\{\,pr_1^2+qr_2^2:\;r_1,r_2\in S\,\}.  \leqno (6.1)
$$
The study of these sets is of interest {\it over both rings and fields.}  
For example, $\,V[\,1,\!1\,]\,$ consists of all sums of two squares in 
$\,S$, and asking if $\,-1\in V[\,1,-d\,]\,$ amounts to solving the 
``negative Pell's equation'' $\,x^2-d\,y^2=-1\,$ over $\,S$.  
If $\,F\,$ is a field of characteristic $\neq 2$, the criterion for the 
splitting of the $F$-quaternion algebra 
$$
\langle \,i,j\;|\;i^2=p,\;j^2=q,\;ij=-ji\,\rangle
\;\;\;\mbox{(where $\,p,\,q \in F\setminus \{0\}$)}  \leqno (6.2)
$$ 
is given by $\,1\in V[\,p,q\,]\,$ (see [La:~p.\,58]). Accordingly, the 
{\it Hilbert symbol\/} $\,(p,q)_F\,$ is defined to be $\,1\,$ or $\,-1$, 
depending on whether or not the quadratic form $\,p\,x^2+q\,y^2\,$ 
represents $\,1\,$ over $\,F$.

\bs
Using the supertrace determinantal formula (4.5), we are able to provide 
in Theorem 6.3 below {\it a matrix-theoretic criterion for a non 
$0$-divisor $\,c\in S\,$ to belong to $\,V[\,p,q\,]\,$ over a commutative 
ring $\,S$.} As far as matrices are concerned, the problem under study 
here is rather distinct from that investigated in the second half of 
[Ta$_2$]. In the latter, Taussky studied the possibility of expressing 
the matrix $\,A\,$ in Theorem 6.3 as a commutator (in the case where 
$\,q=1$), whereas here we are concerned with the {\it factorizations\/} 
of the matrix $\,A\,$ with certain commutator properties.

\sk\nt
{\bf Factorization Theorem 6.3.} {\it  Let $\,A=\begin{pmatrix}0&q\\-p&0
\end{pmatrix} \in R={\mathbb M}_2(S)$, and let $\,c\in S\,$ be a non 
$0$-divisor. The following are equivalent\/}:

\sk\indent
(1) $\;c\in V[\,p,q\,]$.\\
\indent (2) {\it There exist $\,X,\,Y\in R\,$ such that $\,XY=c\cdot A$, 
$\,{\rm det}\,(X)=cp$, $\,{\rm det}\,(Y)=cq$,\\
\indent\;\;\;\;\;\;and $\,{\rm det}\;[X,Y]=-c^2$.}\\
\indent $(3)$ {\it There exist $\,X_1,Y_1\in R\,$ such that 
$\,X_1Y_1=c\cdot A$, $\,{\rm det}\,(X_1)=cq$, $\,{\rm det}\,(Y_1)=cp$,\\
\indent \;\;\;\;\;\;and $\;{\rm det}\;[X_1,Y_1]=-c^2$.}

\bs\nt
{\bf Proof.} $(3)\Rightarrow (1)$. Let $\,X_1,Y_1\in R\,$ be as in $(3)$.
If their supertraces are $\,r\,$ and $\,s$, then 
(since $\,X_1Y_1=c\cdot A\,$ has zero diagonal) Cor.~4.7(1) gives
$\,c^2=-\,{\rm det}\;[X_1,Y_1]=(cq)\,s^2+(cp)\,r^2$.  Cancelling $\,c$,
we get $\,c=pr^2+qs^2\in V[\,p,q\,]$.

\mk\nt
$(1)\Rightarrow (2)$.   (This implication does not require $\,c\,$ to be 
a non $0$-divisor.)  If $\,c\in V[\,p,q\,]$, write $\,c=pr^2+qs^2\,$ for 
some $\,r,s\in S$, and let $\,a=s+pr$, $\,b=r-qs$. Then 
$\,ar-bs=pr^2+qs^2=c$.  For the matrices
$$
X=\begin{pmatrix}a&b\\ps&pr\end{pmatrix}, \;\;\,\mbox{and}\;\;\,
Y=\begin{pmatrix}b&qr\\-a&-qs\end{pmatrix},\leqno (6.4)
$$
we have $\,{\rm det}\,(X)=(ar-bs)\,p=cp$, and 
$\,{\rm det}\,(Y)=(ar-bs)\,q=cq$.  Also,
$$
XY=\begin{pmatrix}0& q\,(ar-bs)\\p\,(bs-ar)&0 \end{pmatrix}
=\begin{pmatrix}0& c\,q\\-c\,p&0 \end{pmatrix}=c\cdot A. \leqno (6.5)
$$
Since $\,{\rm str}\,(X)=a-pr=s\,$ and $\,{\rm str}\,(Y)=b+qs=r$, 
Cor.~(4.7)(1) gives 
$$
{\rm det}\;[X,Y]=-(c\,p)\,r^2-(c\,q)\,s^2=-c^{2}. \leqno (6.6)
$$
\sk\nt
$(2)\Rightarrow (3)$.   (This implication also does not require $\,c\,$ 
to be a non $0$-divisor.)  Given $\,X,Y\,$ as in (2), let 
$\,X_1=Y^{\,\prime}\,$ and $\,Y_1=-X^{\,\prime}\,$ (where the primes
denote the adjoints). Then $\,{\rm det}\,(X_1)={\rm det}\,(Y^{\,\prime})
={\rm det}\,(Y)=cq$, and similarly $\,{\rm det}\,(Y_1)={\rm det}\,(X)=cp$.  
Moreover, $\,X_1Y_1=-Y^{\,\prime}X^{\,\prime}=-(XY)^{\,\prime}=-c\,(-A)=
c\cdot A$.  Finally,
\begin{eqnarray*}
[X_1,Y_1] \!&=&\! [\,Y^{\,\prime},-X^{\,\prime}\,]
=-Y^{\,\prime} X^{\,\prime} + X^{\,\prime}Y^{\,\prime} \\
\!&=&\! -(XY)^{\,\prime} + (YX)^{\,\prime} = -\,[X,Y]^{\,\prime}
\end{eqnarray*}
implies that $\,{\rm det}\;[X_1,Y_1]={\rm det}\;[X,Y]=-c^2$, 
as desired.\qed

\bs\nt
{\bf Remark 6.7.} (A) In the case where both $\,c\,$ and $\,p\,$ are non 
$\,0$-divisors, the condition in (2) that $\,{\rm det}\,(Y)=cq\,$ could 
have been dropped, since it would have followed from $\,XY=c\cdot A\,$
and $\,{\rm det}\,(X)=c\,p$.  However, the present form of the statement 
in (2) is more symmetrical.  (The same remark can be made about the 
statement (3).)

\mk\nt
(B) The implications $(3)\Rightarrow (1)$ and $(2)\Rightarrow (1)$ in 
Thm.~6.3 need not hold if $\,c\in S\,$ is a $0$-divisor.  For instance, 
let $\,S\,$ be the commutative local $\,{\mathbb Q}$-algebra generated 
by $\,x,y\,$ with the relations $\,x^2=y^2=xy=0$, and let $\,c=x$, 
$\,p=q=y$. Then (2) and (3) are trivially satisfied by the matrices 
$\,X=Y=X_1=Y_1=0$. However, $\,c=x\notin V[\,y,y\,]$.

\bs
By further developing the ideas used in the proof of the implication 
$(1)\Rightarrow (2)$ above, we get also the following unexpected 
algebro-geometric result on affine curves over commutative rings.

\bs\nt
{\bf Theorem 6.8.} {\it Given $\,p,q,c\in S$, let $\,C=C_{p,q,c}\,$ be
the plane conic $\,\{(r,s)\in S^2: pr^2+qs^2=c\}$.  Let $\,Q=Q_c\,$ be the
quadric surface $\,\{(x,y,z): xy-z^2=-c^2\}$, and let $\,P=P_{p,q,c}\,$
be the ``vertical plane'' $\,\{(x,y,z):px+qy=-c\}\,$ $($both in $\,S^{\,3})$.
Then there is an affine morphism $\,f: C\rightarrow P\cap Q\,$ defined by}
$$
f(r,s)=\bigl( r\,(2qs-r),\,-s\,(2pr+s),\,rs+pr^2-qs^2 \bigr)\;\;\;\;\,
(\,\forall \,(r,s)\in C\,).  \leqno (6.9)
$$

\mk\nt
{\bf Proof.} Before proceeding with the proof, note that the conic $\,C\,$ 
is nonempty iff $\,c\in V[\,p,q\,]$.  In the case $\,C=\emptyset$, of
course, the statement of the theorem is vacuous. In the following, we may
thus assume that $\,C\neq \emptyset$.  

\mk
Given any point $\,(r,s)\in C\,$ (that is, with $\,pr^2+qs^2=c$), let us 
use the notations and conclusions in the proof of $(1)\Rightarrow (2)$ 
in Thm.~6.3 (recalling that this implication did not require $\,c\,$ to be 
a non $\,0$-divisor in $\,S$).  Since $\,M:=[X,Y]\,$ has trace zero, it can 
be written in the form $\,M=\begin{pmatrix}-z&x\\-y&z\end{pmatrix}$ (for some 
$\,x,y,z\in S$). To compute this matrix, we use the definitions of $\,X,Y\,$ 
in (6.4) (and the fact that $\,XY=c\cdot A\,$)\,:
\begin{eqnarray*}
M=XY-YX &=& \begin{pmatrix}0&qc\\-pc&0\end{pmatrix} 
-\begin{pmatrix} b&qr\\-a&-qs\end{pmatrix}
\begin{pmatrix} a&b\\ps&pr\end{pmatrix} \\
&=& \begin{pmatrix} -(ab+pqrs)&q\,(c-pr^2)-b^2\\a^2-p\,(c-qs^2)&ab+pqrs
\end{pmatrix}.
\end{eqnarray*}
Recalling that $\,a=s+pr\,$ and $\,b=r-qs$, we have
$$
{\left\{ \begin{array}{l}
x=q^2s^2-b^2=(qs+b)\,(qs-b)=r\,(2qs-r),\\
y=p^2r^2-a^2=(pr-a)\,(pr+a)=-s\,(2pr+s),\\
z=ab+pqrs=(s+pr)\,(r-qs)+pqrs=rs+pr^2-qs^2.
\end{array} \right.}
\leqno (6.10)
$$
These are quadratic forms in $\,\{r,s\}\,$ (if we think of 
$\,\{p,q\}\,$ as constants), which  define an affine morphism $\,f\,$
from $\,C\,$ to $\,S^{\,3}$, with the obvious property that 
$\,f(-r,-s)=f(r,s)$.  Furthermore, the fact (from (6.3)(2)) that 
$\,-c^2={\rm det}\,[X,Y] =xy-z^2\,$ implies that $\,f(C)\subseteq Q$. 
Finally, for $\,x,y\in S\,$ as defined above, we have
$$ 
p\,x+q\,y=pr\,(2qs-r)-qs\,(2pr+s)=-(pr^2+qs^2)=-c, \leqno (6.11)
$$
so we have $\,f(C)\subseteq P\,$ also, as desired.\qed
\newpage
\bs\nt
{\bf Remark 6.12.} Some congruence properties of the values of $\,x,y,z\,$
are note-worthy.  For any $\,(r,s)\in C$, (6.10) clearly implies that
$\,x\equiv -r^2\;({\rm mod}\;2q)$, and $\,y\equiv -s^2\;({\rm mod}\;2p)$. 
As for $\,z$, we can rewrite it as follows:
$$
z=rs+(c-qs^2)-qs^2=s\,(r-2qs)+c.  \leqno (6.13)
$$
We did not use this expression for $\,z\,$ in (6.10) since it is not
symmetrical in $\,p\,$ and $\,q\,$ (and also not homogeneous in $\,r\,$
and $\,s$).  However, this new expression does give some additional 
information on $\,z\,$; that is, $\,z\equiv c\;({\rm mod}\;s)$. Similarly, 
we can write $\,z=r\,(s+2pr)-c$, so $\,z\equiv -c\;\,({\rm mod}\;r)\,$
as well.

\bs
To make the meaning of Thm.~6.8 more explicit from the viewpoint of
arithmetic geometry, it is best to work in the case 
$\,S={\mathbb Z}\,$.  In this case, $\,pr^2+qs^2=c\,$ defines a conic 
$\,\overline{C}\subseteq {\mathbb C}^2$, $\,px+qy=-c\,$ defines a 
``vertical plane'' $\,\overline{P} \subseteq {\mathbb C}^3$, while 
$\,xy-z^2=-c^2\,$ defines a quadric surface $\,\overline{Q}\subseteq 
{\mathbb C}^3$.  The map $\,f: \overline{C}\rightarrow\overline{P}\cap 
\overline{Q}\,$ given by the polynomials in (6.10) is then an affine 
morphism {\it defined over $\,{\mathbb Z}$,} taking integer points to 
integer points.  Furthermore, the ring $\,{\mathbb Z}\,$ can be 
replaced throughout by an arbitrary ring of algebraic integers.

\bs\nt
{\bf Example 6.14.} Over $\,S={\mathbb Z}\,$, let $\,p=-3$, $\,q=8$, and 
$\,c=5$.  Obviously, all four points $\,({\pm 1},{\pm 1})\,$ are on the 
conic $\,C$.  Using (6.10), we compute easily that
$$
f(1,1)=f(-1,-1)=(15,5,-10),\;\;\,f(1,-1)=f(-1,1)=(-17,-7,-12), 
\leqno (6.15)
$$
which all lie on the curve $\,P\cap Q$. However, the map $\,f:C\rightarrow
P\cap Q\,$ is not surjective in this example. For instance,
{\it we claim that $\,(15,5,10)\in P\cap Q\,$ is not in $\,f(C)$.} To see 
this, assume for the moment that $\,f(r,s)=(15,5,10)\,$ for some $\,(r,s)
\in C$.  By Remark 6.12, we must have  $\,10\equiv 5\;({\rm mod}\;s)$, so 
$\,s\,$ divides $\,5$. If $\,s={\pm 5}$, then $\,-3r^2+8s^2=5\,$ leads to 
a quick contradiction.  Thus, $\,s={\pm 1}$, and hence also $\,r={\pm 1}$.  
But according to (6.15), the $\,z$-coordinate of $\,f(r,s)\,$ must then be 
either $-10\,$ or $\,-12$, a contradiction.  (In (6.18)(C) below, we'll 
actually give some examples where $\,C=\emptyset$, but $\,P\cap Q\neq 
\emptyset$.) To test the accuracy of the formulas (6.10), it is worthwhile 
to compute a few more image points for the map $\,f$.  For instance, for 
$\,({\pm 3},2)\,$ on the conic $\,C$, we have
$$
f(3,2)=(87,32,-53), \;\;\,\mbox{and}\;\;\, f(-3,2)=(-105,-40,-65),
\leqno (6.16)
$$
which are indeed points in $\,P\cap Q$.  (Recall that the functions 
$\,x,y,z\,$ grow quadratically with respect to the two variables $\,r\,$ 
and $\,s$.)

\bs
We record the following consequence of Theorem 6.8, since we cannot 
locate a reference for it (or for any similar result) in the literature.

\bs\nt
{\bf Corollary 6.17.} {\it For any  $\,p,q\in S$, we have the following\/}:

\mk
(1) {\it If $\,c\in V[\,p,q\,]$, there exist $\,x,y,z\in S\,$ such that 
$\,px+qy=-c\,$ and $\,xy-z^2=-c^2$.} ({\it There also exist $\,x_1,y_1,z_1 
\in S\,$ such that $\,px_1+qy_1=c\,$ and $\,x_1y_1-z_1^2=-c^2$.})\\
\indent (2) {\it If $\;V[\,p,q\,]\,$ contains a unit of $\,S$, then
there exist $\,x_2,y_2,z_2\in S\,$ such that\break
$\,px_2+qy_2=x_2y_2-z_2^2=-1$.}

\bs\nt
{\bf Proof.} (1) follows from Thm.~6.8, and the parenthetical statement 
follows from the main statement by taking $\,(x_1,y_1,z_1)\,$ to be 
$\,(-x,-y,{\pm z})$. 

\mk\nt
(2) Fix a {\it unit\/} $\,c\in V[\,p,q\,]$, and take $\,x,y,z\,$ as 
in (1).  Then $\,x_2=c^{-1}x,\;y_2=c^{-1}y$, and $\,z_2=c^{-1}z\,$ 
satisfy the required conditions.\qed

\bs\nt
{\bf Remark 6.18.} (A) Note that, in (1) above, we cannot say that
``there exist $\,x_3,y_3,z_3\in S\,$ such that $\,px_3+qy_3=c\,$ and 
$\,x_3y_3-z_3^2=c^2$.''  Indeed, for $\,S={\mathbb Z}\,$, take $\,p=-4$, 
$\,q=13$, and $\,c=1$.  We have $\,c\in V[\,p,q\,]\,$ since 
$\,1=-4r^2+13s^2\,$ for $\,(r,s)=(9,5)$.  However, using standard 
software for solving binary quadratic equations (such as [Alp]), 
we can easily check that there {\it do not\/} exist integers 
$\,x_3,y_3,z_3\,$ such that $\,-4x_3+13y_3=x_3y_3-z_{3}^2=1$.


\mk\nt
(B) For $\,S={\mathbb Z}$, the following numerical example shows that,
in case $\,V[\,p,q\,]\,$ contains a unit, say $\,1$, the representation
of $\,1\,$ in the form $pr^2+qs^2$ may involve very large integers 
$\,r\,$ and $\,s$, even though $\,p,q\,$ are pretty small. For instance,
let $\,p=37\,$ and $\,q=-67$. Then $\,1\in V[\,p,q\,]\,$ according to
[Alp], but the smallest solution for $\,37\,r^2-67\,s^2=1\,$ is
$$
r=264,\!638,\!639,\!242,\;\;\;\mbox{and} \;\;\;\;
s=196,\!660,\!308,\!201.  \leqno (6.19)
$$
Confirming our result in Cor.~6.17, [Alp] showed that indeed solutions
exist for the equations $\,37\,x_2-67\,y_2=x_2y_2-z_2^2=-1$.  However, 
the numbers $\,x_2,y_2,z_2\,$ have at least $19$ digits\,!  Of course, 
the specific solution $\,(x,y,z)\,$ constructed from (6.10) by using 
the point $\,(r,s)\,$ in (6.19) is even larger.

\mk\nt
(C) We should also point out that the {\it converse\/} to the main 
statement in (6.17)(1) is not true in general.  For instance, {\it let 
$\,S={\mathbb Z}\,$ again, and take $\,c=1$.} For any $\,p>1\,$ and
$\,q=p+1$, the equations $\,px+qy=xy-z^2=-1\,$ are solved by 
$\,(x,y,z)=(1,-1,0)$, but obviously $\,{\pm 1}\notin V[\,p,p+1\,]$. 
It is, however, possibly more interesting to give an example where 
$\,{\pm 1}\in V[\,p,q\,]\,$ is not simply ruled out ``by absolute 
values''.  For this, we can take, for instance, $\,p=-8,\;q=13$, for 
which the equations $\,px+qy=xy-z^2=-1\,$ are solved by $\,(x,y,z)
=(5,3,4)$. Nevertheless, $\,{\pm 1}\notin V[\,p,q\,]$, in view of the 
fact that $\,{\pm 13}\,$ (or $\,{\pm 5}$) is not a square modulo $\,8$.

\bs
In [KL$_1$], two of the authors study the problem of factorizing a 
matrix $\,A\,$ into a product $\,XY\,$ in such a way that the commutator 
$\,[X,Y]\,$ is invertible. The matrices $\,A\,$ that admit such a 
factorization are said to be {\it reflectable.} For some applications 
of the results (6.3), (6.8), and (6.17) in this section to the study of 
reflectable matrices over commutative rings, see \S5 in [KL$_1$].

\bs

\mk

\bs
\nt Faculty of Mathematics \\
\nt Indian Inst.~of Science Education \& Research, Mohali \\
\nt MGSIPA Transit Campus, Sector 19\\
\nt Chandigarh 160\,019, India

\sk\nt
{\tt dkhurana@iisermohali.ac.in}

\bs
\nt Department of Mathematics \\
\nt University of California \\
\nt Berkeley, CA 94720   

\sk
\nt {\tt lam@math.berkeley.edu}

\bs
\nt Berkeley, CA 94720   

\sk
\nt {\tt shomron@ocf.berkeley.edu}

\bs
\end{document}